\pdfoutput=1

\documentclass[11pt]{amsart}

\usepackage{graphicx, amssymb, amsmath, pstricks, amsthm, pst-node}

\usepackage{pst-all}
\usepackage{latexsym}
\usepackage{amsmath}
\usepackage{epsfig}

\newtheorem{theorem}{Theorem}[section]

\newtheorem{lemma}[theorem]{Lemma}

\newtheorem{corollary}[theorem]{Corollary}

\newtheorem{sublemma}{}[theorem]

\theoremstyle{definition}

\theoremstyle{remark}

\numberwithin{equation}{section}

\newcommand{\ba}{\backslash}

\newcommand{\comment}[1]{}

\begin{document}

\title[Nearly Binary Matroids]{On two classes of nearly binary matroids}

\thanks{The authors were supported by the National Security Agency and by an NSF VIGRE Grant, respectively.}

\author{James Oxley}
\address{Department of
 Mathematics, Louisiana State University, Baton Rouge, Louisiana, USA}
\email{oxley@math.lsu.edu}

\author{Jesse Taylor}
\address{Department of
 Mathematics, Louisiana State University, Baton Rouge, Louisiana, USA}
\email{JTAYL75@math.lsu.edu}

\date{July 24, 2013}

\begin{abstract}
We give an excluded-minor characterization for the class of matroids $M$ in which $M\backslash e$ or $M/e$ is binary for all $e$ in $E(M)$.  This class is closely related to the class of matroids in which every member is binary or can be obtained from a binary matroid by relaxing a circuit-hyperplane.  We also provide an excluded-minor characterization for the second class.
\end{abstract}

\maketitle

\section{Introduction}
\label{introduction}
The class of binary matroids is one of the most widely studied classes of matroids and its members have numerous attractive properties. This motivates the study of classes of matroids whose members are close to being binary. In this paper, we consider one very natural such minor-closed class $\mathcal{Z}$, which consists of those matroids $M$ such that $M\ba e$ or $M/e$ is binary for all elements $e$ of $M$. The main result of the paper is an excluded-minor characterization of $\mathcal{Z}$.  This theorem can be restated in terms of matroid fragility, which has enjoyed a recent surge of research interest. Let $N$ be a matroid. A matroid $M$ is \textit{$N$-fragile} if, for each element $e$ of $E(M)$, at least one of $M\backslash e$ and $M/e$ has no $N$-minor (see, for example, \cite{MWZ}). The class of $N$-fragile matroids is clearly minor-closed. The main result of this paper determines the set of excluded minors for the class of $U_{2,4}$-fragile matroids.  Except where otherwise noted the notation and terminology follow \cite{book}. 

It is well known that if $H$ is a circuit and a hyperplane of a matroid $M$, then there is another matroid $M'$ on $E(M)$ whose bases are the bases of $M$ together with $H$.  We say that $M'$ is obtained from $M$ by \textit{relaxing} the \textit{circuit-hyperplane} $H$ and call $M'$ a \textit{relaxation} of $M$. A class of matroids that arises naturally in determining the excluded minors for $\mathcal{Z}$ is $\mathcal{R}$, those matroids $M$ such that $M$ is binary or $M$ is a relaxation of a binary matroid.

\begin{figure}[ht]
 \begin{center}
  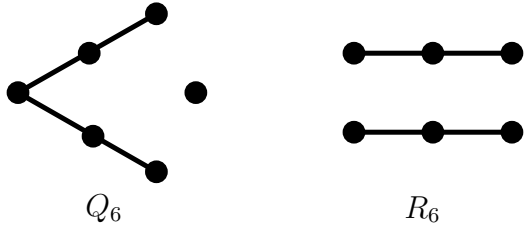
 \end{center}
\caption{Geometric representations of the six-element rank-three matroids $Q_6$ and $R_6$.}
\label{q6r6}
\end{figure}

\begin{figure}[ht]
 \begin{center}
  \includegraphics[scale=.5]{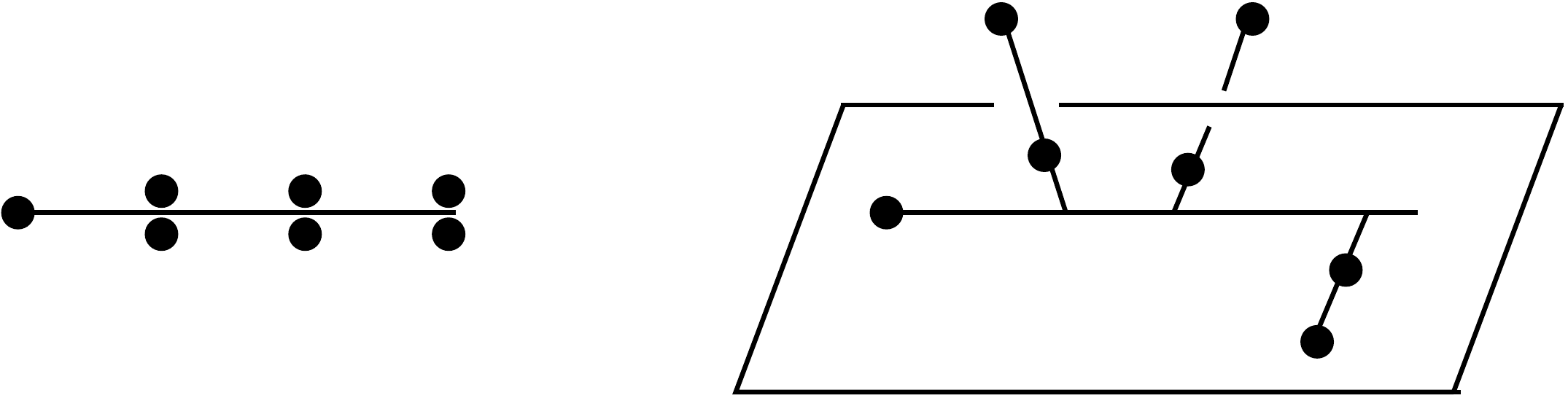}
 \end{center}
\caption{Representations of the matroid $K$ and its (rank-$5$) dual $K^*$.}
\label{Ks}
\end{figure}

The rank-three whirl is denoted by $\mathcal{W}^3$, while $P_6$ is the six-element rank-three matroid that has a single triangle as its only non-spanning circuit.  Let $Q_6$ and $R_6$ be the six-element matroids of rank three for which geometric representations are given in Figure \ref{q6r6}. Evidently $R_6\cong U_{2,4} \oplus_2 U_{2,4}$.  Let $K$ be the seven-element rank-two matroid that is obtained by adding elements in parallel to three of the elements of $U_{2,4}$. The matroid $K$ is depicted with its dual in Figure \ref{Ks}. In Section \ref{prelim}, we note that both $\mathcal{Z}$ and $\mathcal{R}$ are minor-closed and dual-closed classes of matroids and establish some excluded minors of each.  We also introduce Cunningham and Edmonds's canonical tree decomposition of a $2$-connected matroid, along with some preliminaries. 

Let $\mathcal{D}$ denote the collection of all matroids that are obtained from connected binary matroids by relaxing two disjoint circuit-hyperplanes that partition the ground set. The collection $\mathcal{D}$ is in both our sets of excluded minors. Section \ref{main} is devoted to proving the main result and another related result, both of which are stated next.

\begin{theorem}
 The set of excluded minors for the class of matroids $\mathcal{Z}$=$\{M:M\backslash e$ or $M/e$ is binary for all $e$ in $E(M)$\} is $\{Q_6$, $P_6$, $U_{3,6}$, $R_6$, $U_{2,4}\oplus U_{1,1}$, $U_{2,4}\oplus U_{0,1}\}\cup\mathcal{D}$.
\label{nb}
\end{theorem}

\begin{theorem}
 The set of excluded minors for the class $\mathcal{R}$ of matroids $M$ such that $M$ is binary or can be obtained from a binary matroid by relaxing a circuit-hyperplane, is $\{U_{2,5}$, $U_{3,5}$, $K$, $K^*$, $R_6$, $U_{2,4}\oplus U_{1,1}$, $U_{2,4}\oplus U_{0,1}\}\cup\mathcal{D}$.
\label{nb2}
\end{theorem}

For an even integer $r$ exceeding two, let $M_r$ be the rank-$r$ tipless binary spike, that is, the vector matroid of the binary matrix $[I_r|J_r-I_r]$ where $J_r$ is the matrix of all ones. Labeling the columns of this matrix $e_1,e_2,\ldots, e_{2r}$ in order, we see that $\{e_2,e_3,\ldots,e_r,e_{r+1}\}$ and its complement are both circuit-hyperplanes of $M_r$.  By relaxing these circuit-hyperplanes, we obtain a member of $\mathcal{D}$.  Thus the sets of excluded minors in Theorems \ref{nb} and \ref{nb2} are both infinite. However, these doubly relaxed spikes are not the only members of $\mathcal{D}$.  In Section \ref{compD}, we further discuss the complexity of $\mathcal{D}$.

As $\mathcal{D}$ shows, the class of matroids that can be obtained from binary matroids by relaxing at most two circuit-hyperplanes does contain an infinite antichain. Geelen, Gerards, and Whittle announced in 2009 that the class of binary matroids itself contains no infinite antichains. These observations raise the interesting question, which was asked by a referee of this paper, as to whether or not the class $\mathcal{R}$ contains an infinite antichain. It is not difficult to check using, for example, \cite[Lemma 2.6]{hine}, that $\mathcal{Z}$ contains an infinite antichain if and only if $\mathcal{R}$ does.

\section{Preliminaries}
\label{prelim}

This section first notes that both $\mathcal{Z}$ and $\mathcal{R}$ are minor- and dual-closed, and then determines some excluded minors for each class.  

\begin{lemma}
The classes $\mathcal{Z}$ and $\mathcal{R}$ are both closed under duality and the taking of minors.


\end{lemma}

This lemma is immediate for $\mathcal{Z}$ and is a straightforward consequence of the following result of Kahn \cite{kahn} for $\mathcal{R}$ (see also \cite[p. 115]{book}).

\begin{lemma}
Let $X$ be a circuit-hyperplane of a matroid $M$ and let $M'$ be the matroid obtained from $M$ by relaxing $X$.  Then $(M')^*$ is obtained from $M^*$ by relaxing the circuit-hyperplane $E(M)-X$ of the latter.  Moreover, when $e \in E(M)-X$, $M/e$ and $M'/e$ are equal and, unless $M$ has $e$ as a coloop, $M'\backslash e$ is obtained from $M\ba e$ by relaxing the circuit-hyperplane $X$ of the latter, and the dual situation holds when $e\in X$.
\label{relax}
\end{lemma}

It is not difficult to deduce from the above result that the class $\mathcal{R}$ is contained in the class $\mathcal{Z}$.  We say a matroid $N$ is a \textit{series extension} of a matroid $M$ if $M=N/T$ and every element of $T$ is in series with some element of $M$.  We call $N$ a \textit{parallel extension} of $M$ if $N^*$ is a series extension of $M^*$.  Note that this differs from the terminology used in \cite{book}.  The following result from \cite{list} will be used extensively throughout the paper.

\begin{theorem}
A matroid $M$ is non-binary and in $\mathcal{Z}$ if and only if
	\begin{enumerate}
		\item [(i)] both $r(M)$ and $r^*(M)$ exceed two and $M$ can be obtained from a connected binary matroid by relaxing a circuit-hyperplane; or
		\item [(ii)] $M$ is isomorphic to a parallel extension of $U_{2,n}$ for some $n$ $\geq$ 5; or
		\item [(iii)] $M$ is isomorphic to a series extension of $U_{n-2,n}$ for some $n$ $\geq$ 5; or 
		\item [(iv)] $M$ can be obtained from $U_{2,4}$ by series extension of a subset $S$ of $E(U_{2,4})$ and parallel extension of a disjoint subset $T$ of $E(U_{2,4})$ where $S$ or $T$ may be empty.
	\end{enumerate}
\label{class}
\end{theorem}

Let $EX(\mathcal{M})$ denote the class of excluded minors for a class of matroids $\mathcal{M}$. Some excluded minors for $\mathcal{Z}$ and $\mathcal{R}$ are easy to identify.  We omit the routine argument that establishes the following.

\begin{lemma}
 The matroids $U_{2,4}\oplus U_{1,1}$, $U_{2,4}\oplus U_{0,1}$, and $R_6$ are in both $EX(\mathcal{Z})$ and $EX(\mathcal{R})$.  
\label{share}
\end{lemma}

The following three results will also be useful, the first is from \cite{kahn}; the second is elementary; the third follows from the first two.

\begin{lemma}
Let $M'$ be obtained from $M$ by relaxing a circuit-hyperplane.
\begin{enumerate}
 \item [(i)] If $M$ is connected, then $M'$ is non-binary; and
 \item [(ii)] if $M$ is $n$-connected, then so is $M'$.
\end{enumerate}
\label{kahnconn}
\end{lemma}

\begin{lemma}
 The only disconnected matroids having a circuit-hyperplane are $U_{n-1,n}\oplus U_{1,k}$, for integers $n$, $k$ $\geq 1$.
\label{nonbin}
\end{lemma}

\begin{corollary}
 Let $M$ be a binary matroid, $H$ be a circuit-hyperplane of $M$, and $M'$ be obtained from $M$ by relaxing $H$. Then $M'$ is binary if and only if $M$ is $U_{n-1,n}\oplus U_{1,k}$, for integers $n$, $k$ $\geq 1$.
\label{oxwhit}
\end{corollary}

Note that, in Lemma \ref{nonbin} and Corollary \ref{oxwhit}, the disconnected matroids are graphic and carry the name \textit{enlarged $1$-wheels} in \cite{whit}.

Recall, $\mathcal{D}$ is the collection of all matroids that are obtained from connected binary matroids by relaxing two disjoint circuit-hyperplanes that partition the ground set.

\begin{lemma}
All matroids in $\mathcal{D}$ are in both $EX(\mathcal{Z})$ and $EX(\mathcal{R})$.
\begin{proof}
Take a matroid $M_2$ in $\mathcal{D}$.  Let $X$ and $Y$ be the disjoint circuit-hyperplanes of the connected binary matroid $M$ that are relaxed to obtain $M_2$. Let $M_X$ and $M_Y$ denote the matroids obtained from $M$ by relaxing $X$ and $Y$, respectively, and take $e$ in $E(M_2)$.  Note that the case with $e\in X$ is symmetric to the case with $e\in Y$; both $\mathcal{Z}$ and $\mathcal{R}$ are dual-closed classes, and since $X$ and $Y$ are complementary circuit-hyperplanes of $M$, they are so for $M^*$ as well.

Suppose $e\in X$.  By Lemma \ref{relax}, $M_2/e$ is obtained from $M_Y/e$ by relaxing the circuit-hyperplane $X-e$ of the latter and $M_Y/e$=$M/e$.  If $M/e$ is connected, then $M_2/e$ is non-binary by Lemma \ref{kahnconn}.  Now assume $M/e$ is disconnected.  Then $M/e$=$U_{n-1,n}\oplus U_{1,k}$ for some $n$, $k$ $\geq1$, by Lemma \ref{nonbin}.  But $Y$ is a spanning circuit in $M/e$, which is a contradiction since $M/e$ has no spanning circuits. We conclude that $M_2/e$ is non-binary.  By symmetry and duality the same argument holds for $M_2\ba e$, and for both $M_2/f$ and $M_2\ba f$ when $f\in Y$.

Any deletion $M_2\ba z$ equals $M_Y\ba z$ or $M_X\ba z$. By symmetry we only need to consider the case with $z\in X$. The matroid $M_Y\ba z$ can be obtained by relaxing a circuit-hyperplane in a binary matroid. By duality, the same holds for $M_2/z$. Therefore any minor of $M_2$ is in $\mathcal{R}$ and so is in $\mathcal{Z}$.  Thus $M_2$ is in $EX(\mathcal{R})$ and in $EX(\mathcal{Z})$.
\end{proof}
\label{D} 
\end{lemma}

The next two lemmas list matroids that are excluded minors for exactly one of $\mathcal{R}$ and $\mathcal{Z}$.  Their routine proofs are omitted.

\begin{lemma}
The matroids $U_{2,5}$, $U_{3,5}$, $K$, and $K^*$ are excluded minors for the class $\mathcal{R}$.

\label{R}
\end{lemma}

\begin{lemma}
 The matroids $Q_6$, $P_6$, and $U_{3,6}$ are excluded minors for the class $\mathcal{Z}$.
\label{Z}
\end{lemma}


A class $\mathcal{N}$ of matroids is \textit{$1$-rounded} \cite{2round} if every member of $\mathcal{N}$ is connected and, whenever $e$ is an element of a connected matroid $M$ having an $\mathcal{N}$-minor, $M$ has an $\mathcal{N}$-minor using $e$.  The following three results will be useful in our proofs, they come from \cite{bixby}, \cite{2round}, and \cite{round}, respectively.

\begin{lemma}
 The set $\{U_{2,4}\}$ is $1$-rounded.
\label{bix}
\end{lemma}

\begin{lemma}
 The set $\{M(K_4)$, $U_{2,4}\}$ is $1$-rounded.
\label{roun}
\end{lemma}

\begin{lemma}
 The set $\{\mathcal{W}^3$, $P_6$, $Q_6$, $U_{3,6}\}$ is $1$-rounded.
\label{again}
\end{lemma}

Next we introduce Cunningham and Edmonds's tree decomposition for connected matroids \cite{cun}.  Our treatment of this material follows \cite[pp. 307--310]{book}.  A \textit{matroid-labeled tree} is a tree $T$ with vertex set $\{M_1$, $M_2$, $\ldots$, $M_k\}$ for some positive integer $k$ such that
\begin{enumerate}
 \item [(i)] each $M_i$ is a matroid;
 \item [(ii)] if $M_{j_1}$ and $M_{j_2}$ are joined by an edge $e_i$ of $T$, then $E(M_{j_1}) \cap E(M_{j_2})$=$\{e_i\}$, and $\{e_i\}$ is not a separator of $M_{j_1}$ or $M_{j_2}$; and
 \item [(iii)] if $M_{j_1}$ and $M_{j_2}$ are non-adjacent, then $E(M_{j_1})\cap E(M_{j_2})$ is empty.
\end{enumerate}

Let $e$ be an edge of a matroid-labeled tree $T$ and suppose $e$ joins vertices labeled by $M_1$ and $M_2$.  Suppose that we contract $e$ and relabel by $M_1 \oplus_2 M_2$ the composite vertex that results by identifying the endpoints of $e$.  Then we retain a matroid-labeled tree and we denote this tree by $T/e$.  This process can be repeated and since the operation of $2$-sum is associative, for every subset $\{e_{i_1}, e_{i_2}, \ldots, e_{i_m}\}$ of $E(T)$, the matroid-labeled tree $T/e_{i_1}, e_{i_2}, \ldots, e_{i_m}$ is well-defined.

A \textit{tree decomposition} of a $2$-connected matroid $M$ is a matroid-labeled tree $T$ such that if $V(T)$=$\{M_1$, $M_2$, $\ldots$, $M_k\}$ and $E(T)$=$\{e_1$, $e_2$, $\ldots$, $e_{k-1}\}$, then
\begin{enumerate}
 \item [(i)] $E(M)$=$(E(M_1)\cup E(M_2)\cup \cdots \cup E(M_k))-\{e_1,e_2,\ldots,e_{k-1}\}$;
 \item [(ii)] $E(M_i)\geq 3$ for all $i$ unless $|E(M)| < 3$, in which case $k$=1 and $M_1$=$M$; and
 \item [(iii)] $M$ labels the single vertex of $T/e_1,e_2,\ldots,e_{k-1}$.
\end{enumerate}

In general, a tree decomposition of a matroid is not unique.  However, Cunningham and Edmonds were able to guarantee uniqueness of the \textit{canonical tree decomposition} described in the following theorem from \cite{cun}. 

\begin{theorem}
 Each $2$-connected matroid $M$ has a tree decomposition $T$ in which every vertex is labeled by a $3$-connected matroid, $U_{m-1,m}$ for some $m\geq 3$, or $U_{1,n}$ for some $n\geq 3$. Moreover, there are no two adjacent vertices that are both labeled by uniform matroids of rank one or are both labeled by uniform matroids of corank one, and $T$ is unique to within a relabeling of its edges.
\end{theorem}

The canonical tree decomposition provides a unique way to break up a $2$-connected matroid $M$ into $3$-connected pieces, uniform matroids of rank one, and uniform matroids of corank one.  Moreover, we can reconstruct $M$ from these pieces using the $2$-sum operation with the common elements between matroids as basepoints.  A basic property of the $2$-sum operation is that $M_1$ and $M_2$ are minors of $M_1 \oplus_2 M_2$.  The following result is well known; its routine proof is omitted.

\begin{lemma}
 Let $M_1$ and $M_2$ label vertices in a tree decomposition $T$ of a connected matroid $M$.  Let $P$ be the path in $T$ joining $M_1$ and $M_2$, and let $p_1$ and $p_2$ be the edges of $P$ meeting $M_1$ and $M_2$ respectively.  In other words, $p_1$ and $p_2$ are basepoints for $2$-sums in the reconstruction of $M$.  Then $M$ has a minor isomorphic to the $2$-sum of $M_1$ and $M_2$, where $p_1$=$p_2$ is the basepoint of the $2$-sum.
\label{shrink}
\end{lemma}

The following two results will also be needed. The first is basic and its proof is omitted.  The second result comes from \cite{min3}.

\begin{lemma}
 The class of binary matroids is closed under the operation of $2$-sum.
\label{bin2sum}
\end{lemma}

\begin{lemma}
 The following statements are equivalent for a $3$-connected matroid $M$ having rank and corank at least three:
\begin{enumerate}
 \item [(i)] $M$ has a $U_{2,5}$-minor;
 \item [(ii)] $M$ has a $U_{3,5}$-minor;
 \item [(iii)] $M$ has a minor isomorphic to one of $P_6$, $Q_6$, or $U_{3,6}$.
\end{enumerate}
\label{extra}
\end{lemma}

\section{Main Result}
\label{main}

In this section we prove the main results of the paper, Theorems \ref{nb} and \ref{nb2}.  We begin by finding all the disconnected excluded minors of each class. Due to the similarity of the proofs for each class, we combine the arguments where possible.

\begin{lemma}
 Suppose $\mathcal{U}$ $\in$ $\{\mathcal{Z}, \mathcal{R}\}$.  The only disconnected members of $EX(\mathcal{U})$ are $U_{2,4} \oplus U_{1,1}$ and $U_{2,4} \oplus U_{0,1}$.
\begin{proof}
 By Lemma \ref{share}, both matroids are in $EX(\mathcal{U})$.  Now let $M$ be an arbitrary disconnected member of $EX(\mathcal{U})$.  As $M$ is non-binary and disconnected, it has distinct components $M_1$ and $M_2$ where $M_1$ is non-binary.  Since $M_1$ has a $U_{2,4}$-minor and $M_2$ has a $U_{0,1}$- or $U_{1,1}$-minor, the lemma follows.  
\end{proof}
\label{discon}
\end{lemma}

The following result from \cite{min3} will be useful in our proofs.

\begin{theorem}
Let $M$ be a $3$-connected matroid having rank and corank exceeding two.  
  \begin{enumerate}
   \item [(i)] If $M$ is binary, then $M$ has an $M(K_4)$-minor.
   \item [(ii)] If $M$ is non-binary, then $M$ has one of $\mathcal{W}^3$, $Q_6$, $P_6$, and $U_{3,6}$ as a minor.
  \end{enumerate}
\label{3connminor}
\end{theorem}

Before finding the complete list of $2$-connected excluded minors, we need the following lemmas.  The first lemma comes from \cite[Section 1.5, Exercise 14]{book}; its proof is omitted.

\begin{lemma}
 The following statements are equivalent for a matroid $M$:
\begin{enumerate}
 \item [(a)] $M$ is a relaxation of some matroid,
 \item [(b)] $M$ has a basis $B$ such that $B\cup e$ is a circuit of $M$ for every $e$ in $E(M)-B$ and neither $B$ nor $E(M)-B$ is empty.
\end{enumerate}

\label{relaxbase}
\end{lemma}

\begin{lemma}
 Let $M$ be a matroid that can be obtained from a binary matroid $N$ by relaxing a circuit-hyperplane $X$ of the latter.  If $M$ contains a $\mathcal{W}^k$-minor for some $k\geq 3$, then, in every $\mathcal{W}^k$-minor of $M$, the rim elements are contained in $X$ and no element of $X$ is a spoke.
\begin{proof}
 Let $M_1$ be a $\mathcal{W}^k$-minor of $M$.  If $e$ is in the rim of $M_1$, then $M_1/e$ is non-binary.  But, for all $f$ in $E(M)-X$, by Lemma \ref{relax}, $M/f$ is binary.  Therefore $e\in X$.  The assertion about spokes follows by duality.
\end{proof}
\label{bases}
\end{lemma}

\begin{lemma}
 Let $M$ be a connected non-binary matroid.  Either $M$ has an $R_6$-, $U_{2,4}\oplus U_{0,1}$-, or $U_{2,4}\oplus U_{1,1}$-minor, or $M$ is obtained from a $3$-connected non-binary matroid $M_0$ by parallel and series extension of disjoint subsets $T$ and $S$ of $E(M_0)$, where both $S$ and $T$ are possibly empty.
\begin{proof}
Consider the canonical tree decomposition $T$ of $M$. As $M$ is non-binary, by Lemma \ref{bin2sum} there must be a non-binary matroid $M_0$ in $T$. Assume there is another vertex labeled by a non-binary matroid $M_1$.  Then, by Lemma \ref{shrink}, we see that $M$ has an $M_0\oplus_2 M_1$-minor.  Let $p_1$ be the basepoint of this $2$-sum.  Each of $M_0$ and $M_1$ is connected and non-binary, so by Lemma \ref{bix} each of $M_0$ and $M_1$ has a $U_{2,4}$-minor that uses $p_1$.  Thus $M$ has an $R_6$-minor, and the lemma holds when $M_1$ exists.  

We may now assume that $M_0$ is the unique non-binary matroid labeling a vertex of $T$.  Suppose there is a vertex labeled by a $3$-connected binary matroid $M_2$ with at least four elements.  Then $M$ has an $M_0 \oplus_2 M_2$-minor.  Now $M_0$ has a $U_{2,4}$-minor and, as $M_2$ is $3$-connected and binary, Theorem \ref{3connminor} tells us that $M_2$ has an $M(K_4)$-minor.  Let $p_2$ be the basepoint of $M_0 \oplus_2 M_2$.  As above, $M_0$ has a $U_{2,4}$-minor using $p_2$.  By Lemma \ref{roun}, $M_2$ has an $M(K_4)$-minor using $p_2$.  Thus $M$ has a $U_{2,4} \oplus_2 M(K_4)$-minor and therefore has a $U_{2,4} \oplus U_{1,1}$-minor.  Hence the lemma holds when $M_2$ exists.

We may now assume all matroids other than $M_0$ labeling vertices in $T$ are $U_{1,n}$ or $U_{m-1,m}$ for varying $n$, $m$ $\geq 3$.  If we have a path in $T$ beginning at $M_0$ that has the form $M_0$---$U_{m-1,m}$---$U_{1,n}$, then $M$ has a $U_{2,4} \oplus U_{0,1}$-minor.  By duality, we may not have a path of the form $M_0$---$U_{1,n}$---$U_{m-1,m}$.  Therefore we may assume the only non-trivial paths beginning at $M_0$ in $T$ are of the form $M_0$---$U_{m-1,m}$, or $M_0$---$U_{1,n}$.  In other words, $M$ is obtained from $M_0$ by parallel and series extension of disjoint subsets of $E(M_0)$.
\end{proof}
\label{decomp}
\end{lemma}

Recall that the matroid $K$ is the matroid obtained from $U_{2,4}$ by adding elements in parallel to three of its elements.

\begin{lemma}
The matroid $R_6$ is the only connected, but not $3$-connected, member of $EX(\mathcal{Z})$.  The connected, but not $3$-connected, members of $EX(\mathcal{R})$ are $R_6$, $K$, and $K^*$.
\begin{proof}
Suppose $\mathcal{U}$ $\in$ $\{\mathcal{Z},\mathcal{R}\}$.  By Lemmas \ref{share} and \ref{R}, $R_6$ is in $EX(\mathcal{U})$ and $K$ and $K^*$ are in $EX(\mathcal{R})$.  Let $M$ be a $2$-connected member of $EX(\mathcal{U})$ that is not $3$-connected and is not $R_6$, $K$, or $K^*$.  By Lemma \ref{decomp}, $M$ is obtained from a $3$-connected non-binary matroid $M_0$ by parallel and series extension of disjoint subsets $T$ and $S$ of $E(M_0)$ where $S\cup T \neq \emptyset$. 

Let $M_0$ $\cong$ $U_{2,4}$.  If $\mathcal{U}$=$\mathcal{Z}$, then $M$ is in $\mathcal{U}$, as it satisfies (iv) in Theorem \ref{class}, which is a contradiction, so let $\mathcal{U}$=$\mathcal{R}$. As $M$ has neither $K$ nor $K^*$ as a minor, both $S$ and $T$ have size less than three.  By duality, we may assume that $0\leq|T|\leq|S|\leq 2$. In each case, $M$ can be realized as a relaxation of a binary matroid. For example, when $|S|=|T|=2$, assume the non-trivial series classes have sizes $s_1$ and $s_2$, and the non-trivial parallel classes have sizes $p_1$ and $p_2$.  We can obtain $M$ by relaxing the circuit-hyperplane in $M(G)$ where $G$ is a graph on three vertices $\{a, b, c\}$ with $p_1$ parallel edges between $a$ and $c$, $p_2$ parallel edges between $b$ and $c$, and two internal vertex disjoint paths with sizes $s_1$ and $s_2$ between $a$ and $b$.  The other cases can be checked similarly.  We deduce contradictorily that $M\in \mathcal{U}$.

We may now assume $|E(M_0)|$ $\geq$ 5 and consider $\mathcal{U}$ $\in$ $\{\mathcal{Z},\mathcal{R}\}$.  By switching to the dual if necessary, we may also assume that $M$ has at least one non-trivial parallel class and let $\{x,y\}$ be in that class.

\begin{sublemma}
The matroid $M\ba x$ can be obtained from a binary matroid by relaxing a circuit-hyperplane.
\begin{proof}
We know $M\ba x$ is in $\mathcal{U}$.  Thus it satisfies one of (i)-(iv) in Theorem \ref{class}.  If $M\ba x$ satisfies (i), then the result follows. Assume $M\ba x$ satisfies (ii).  Then $M\ba x$ is a parallel extension of $U_{2,n}$, for some $n$ $\geq$ 5.  Hence $M$ is also a parallel extension of this matroid and $M\in \mathcal{Z}$, and $M$ has a $U_{2,5}$-minor, which contradicts Lemma \ref{R} if $\mathcal{U}$=$\mathcal{R}$. Next assume that $M\ba x$ satisfies (iii).  Then $M\ba x$ is a series extension of $U_{n-2,n}$ for some $n$ $\geq$ 5, and $M$ is a parallel extension of this series extension.  Then $M\ba x$, and hence $M$, contains the excluded minor $U_{2,4} \oplus U_{0,1}$.  Lastly, assume $M\ba x$ satisfies (iv).  Let $U=E(U_{2,4})-S-T$, where $S$ and $T$ are as defined in Theorem \ref{class}. Recall that $\{x,y\}$ is a circuit of $M$.  If $y$ is in a non-trivial series class of $M\ba x$, then $M$ contains the excluded minor $U_{2,4}\oplus U_{0,1}$, so $y\in T\cup U$.  Therefore $M$ satisfies (iv), so $M\in \mathcal{Z}$ and we assume $\mathcal{U}$=$\mathcal{R}$.  As $M$ has neither $K$ nor $K^*$ as a minor, $|S|<3$ and $|T\cup y|<3$.  As noted above, in these cases $M$ can be realized as a relaxation of a binary matroid. 
\end{proof}
\label{sublem}
\end{sublemma}
%

If $r(M_0)=2$, then $M_0$ has a $U_{2,5}$-minor, so assume $\mathcal{U}$=$\mathcal{Z}$.  It is not hard to check that we get a contradiction in this case by establishing that either $M\in \mathcal{Z}$, or $M$ contains a $U_{2,4} \oplus U_{1,1}$-minor.  Thus we may assume $\mathcal{U}$ $\in$ $\{\mathcal{Z},\mathcal{R}\}$, $r(M_0)\geq 3$ and, by duality, $r^*(M_0)\geq 3$.  As $M_0$ is non-binary and $3$-connected, and all of $P_6$, $Q_6$, and $U_{3,6}$ are either in $EX(\mathcal{U})$ or contain members of $EX(\mathcal{U})$, Theorem \ref{3connminor} implies that $M_0$ contains a $\mathcal{W}^3$-minor.  By \ref{sublem}, $M\ba x$ is a relaxation of a binary matroid $N$, so let $B$ be the circuit-hyperplane relaxed in $N$ to produce $M\ba x$.  Assume $y\notin B$ and let $N_1$ be obtained from $N$ by adding $x$ back in parallel to $y$.  Then $B$ is a circuit-hyperplane of $N_1$ whose relaxation is $M$, a contradiction.

We may now assume $y$ $\in$ $B$. Since $M_0$ has a $\mathcal{W}^3$-minor and no $P_6$-, $Q_6$-, or $U_{3,6}$-minor, by Lemma \ref{again} $M_0$ has a $\mathcal{W}^3$-minor $M_y$ using $y$.  By Lemma \ref{bases}, we know that $y$ is a rim element of $M_y$.  This implies that $M$ has a $\mathcal{W}^3$-minor in which one of the rim elements is replaced by the parallel class containing $\{x,y\}$. This is a contradiction since it implies $M$ has a $U_{2,4} \oplus U_{0,1}$-minor.
\end{proof}
\label{conn}
\end{lemma}
 
In finding the complete list of $3$-connected excluded minors for each class, we use the following lemma.

\begin{lemma}
If a matroid $N'$ is obtained from a non-binary matroid $N$ by relaxing a circuit-hyperplane $X$, then 
\begin{enumerate}
 \item [(i)] $N'$ has a $U_{2,5}$- or $U_{3,5}$-minor; or
 \item [(ii)] $N'$ has a matroid in the class $\mathcal{D}$ as a minor.
\end{enumerate}
\begin{proof}
As $N$ is non-binary and has a circuit-hyperplane, $|E(N)|\geq 5$.  If $|E(N)|=5$, then either $N$ is $U_{2,4} \oplus_2 U_{1,3}$, in which case $N'$ is $U_{2,5}$, or $N$ is $U_{2,4} \oplus_2 U_{2,3}$, in which case $N'$ is $U_{3,5}$.  Thus the result holds if $|E(N)|=5$.  Now assume that the result holds for $|E(N)|<k$, and consider the case where $|E(N)|=k\geq 6$.  If $r(N)=2$, then $N'$ has a $U_{2,5}$-minor and the result holds. Dually, the result holds if $r^*(N)=2$, so assume $r(N),r^*(N)\geq 3$.

Take $e\in X$ and consider $N/e$.  By Lemma \ref{relax}, $N'/e$ is obtained from $N/e$ by relaxing $X-e$.  If $N/e$ is non-binary, then we invoke the induction hypothesis to see that the result holds.  Hence $N/e$ is binary for all $e\in X$.  By duality, $N\ba e$ is binary for all $e\not \in X$.  Thus, for every $e\in E(N)$, at least one of $N\ba e$ and $N/e$ is binary.  By Theorem \ref{class}, we deduce that one of (i)-(iv) holds for $N$.

As $r(N), r^*(N)\geq 3$, we know $N$ cannot satisfy (ii) or (iii).  Assume $N$ satisfies (iv).  If $|S|=0$ or $|T|=0$, we contradict our rank or corank assumptions, so $|S|, |T|\geq 1$.  It is straightforward to check that $N$ cannot have a circuit-hyperplane, which is a contradiction.

Finally, assume $N$ satisfies (i).  Then $N$ can be obtained from some connected binary matroid $M$ by relaxing a circuit-hyperplane $Y$ in $M$.  Assume $X\cap Y\neq\emptyset$ and take $e\in X\cap Y$.  Then $N/e$ is binary and is obtained from the binary matroid $M/e$ via relaxation.  By Corollary \ref{oxwhit}, $M/e\cong U_{n-1,n}\oplus U_{1,k}$ for some $n, k\geq 1$, and $N/e \cong U_{n,n+1}\oplus_2 U_{1,k+1}$.  However, this implies $N/e$ has no circuit-hyperplane unless $n=2$ and $k=2$, so we assume these values for $n$ and $k$.  But this means $N'/e\cong U_{2,4}$, which is a contradiction since $r(N'), r^*(N')\geq 3$.  Thus $X\cap Y=\emptyset$ and, by duality, $(E(M)-X)\cap(E(N)-Y)=\emptyset$.  As both $(X, E(N)-X)$ and $(Y, E(N)-Y)$ partition the ground set, $X=E(N)-Y$ and $E(N)-X=Y$.  Hence $N'$ is obtained from the connected binary matroid $M$ by relaxing the two disjoint circuit-hyperplanes $X$ and $Y$, so $N'$ is in $\mathcal{D}$ and the result holds.
\end{proof}
\label{import}
\end{lemma}

\begin{figure}[ht]
 \begin{center}
  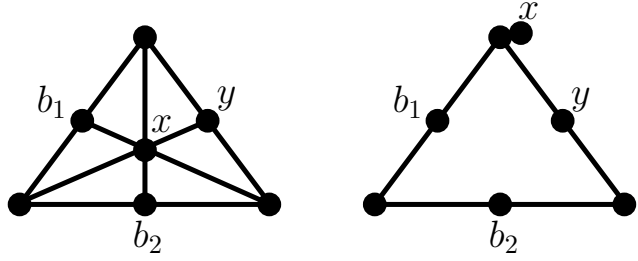
 \end{center}
\caption{Geometric representations of $F^-_7$, and a parallel extension of a spoke element of $\mathcal{W}^3$.}
\label{ext}
\end{figure}

\begin{lemma}
The complete list of $3$-connected members of $EX(\mathcal{Z})$ is $Q_6$, $P_6$, $U_{3,6}$, and the matroids in $\mathcal{D}$.  The complete list of $3$-connected members of $EX(\mathcal{R})$ is $U_{2,5}$, $U_{3,5}$, and the matroids in $\mathcal{D}$
\begin{proof}
Suppose $\mathcal{U}$ $\in$ $\{\mathcal{Z},\mathcal{R}\}$.  Let $M$ be a $3$-connected excluded minor of $\mathcal{U}$ that is not $Q_6$, $P_6$, $U_{3,6}$, $U_{2,5}$, $U_{3,5}$, or any of the matroids in $\mathcal{D}$.  Clearly $r(M)\geq 3$ and $r^*(M)\geq 3$.  Either (a) $M$ is a relaxation of a non-binary matroid; or (b) $M$ is not a relaxation of any matroid at all. Case (a) follows immediately by Lemmas \ref{import} and \ref{extra}.

Now consider case (b). By Theorem \ref{3connminor}, $M$ must contain one of $\mathcal{W}^3$, $Q_6$, $P_6$, and $U_{3,6}$. As all of these except $\mathcal{W}^3$ contain excluded minors of $\mathcal{U}$, we know that $M$ has a $\mathcal{W}^3$-minor.  Let $\mathcal{W}^k$ be the largest whirl-minor of $M$.  We use Seymour's Splitter Theorem \cite{seymour} to grow $M$ from $\mathcal{W}^k$.  Let $x$ be the element added with the last move.  By duality, we may assume that $x$ is added via extension.  Thus $M\ba x$ is a non-binary $3$-connected member of $\mathcal{U}$.  If $\mathcal{U}$=$\mathcal{R}$, then $M\ba x$ is a relaxation of a binary matroid. If $\mathcal{U}$=$\mathcal{Z}$, then $M\ba x$ satisfies one of (i)-(iv) in Theorem \ref{class}.  As $M\ba x$ is $3$-connected, it cannot satisfy (ii)-(iv). Hence, in both cases, $M\ba x$ is a relaxation of a binary matroid $N_1$.

Let $B$ be the special basis in $M\ba x$ that is a circuit-hyperplane in $N_1$.  For all $e \in$ $E(M\ba x)-B$, the set $B \cup e$ is a circuit in $M\ba x$.  Now $B$ is also a basis of $M$ and $B\cup e$ is a circuit of $M$ for all $e \in$ $E(M)-(B\cup x)$.  If $B\cup x$ is a circuit of $M$, then $M$ can be realized as a relaxation of some matroid by Lemma \ref{relaxbase}, which is a contradiction.  Thus there is some $y \in B$ such that $y$ is not in the circuit contained in $B\cup x$.  Now, by Lemma \ref{relax}, $M\ba x/y$ can be obtained by relaxing the circuit-hyperplane $B-y$ in $N_1/y$.

Assume that $M\ba x/y$ is binary.  It follows from Corollary \ref{oxwhit} that $M\ba x/y$ can be obtained from a circuit $C$ by adding some, possibly empty, set of elements in parallel with some element $z$ of $C$ where $C-z=B-y$.  As $M\ba x$ is $3$-connected, it has no non-trivial series classes.  Hence $|C-z|=1$, so $r(M\ba x/y)=1$, which contradicts the fact that $r(M)\geq 3$.  Therefore $M\ba x/y$ must be non-binary, and so $M/y$ is also non-binary.

\begin{sublemma}
The matroid $M/y$ can be obtained from a binary matroid via relaxation.  
\begin{proof}
This is certainly true if $\mathcal{U}$=$\mathcal{R}$, so assume $\mathcal{U}$=$\mathcal{Z}$.  Then $M/y$ satisfies one of (i)-(iv) in Theorem \ref{class}. First note that $M/y$ cannot satisfy (iii), because a connected single-element coextension of a series extension of $U_{n-2,n}$ has corank two, and so has no $\mathcal{W}^3$-minor.  Assume $M/y$ satisfies (ii).  Then $r(M)=3$. As $M$ has a $\mathcal{W}^3$-minor, it is not hard to check that we must coextend $M/y$ by $y$ in a way that creates a matroid having a $U_{2,4}\oplus U_{1,1}$-, $Q_6$-, or $P_6$-minor.  Now assume $M/y$ satisfies (iv).  As $M$ is $3$-connected, $M/y$ cannot have any non-trivial series classes.  A routine check shows that either $M$ contains an excluded minor or $M$ can be realized as a relaxation of a binary matroid, both of which are contradictions.  Thus (i) holds, and so does the result.
\end{proof}
\label{newbie} 
\end{sublemma}

We now revert to working in generality, where $\mathcal{U}$ $\in$ $\{\mathcal{Z}$,$\mathcal{R}\}$.  We know that $M\ba x$ is obtained by relaxing a circuit-hyperplane $B$ in a binary matroid $N_1$, and $M/y$ is obtained by relaxing a circuit-hyperplane $B'$ in a binary matroid $N_2$.   We show next that 

\begin{sublemma}
$B'$=$B-y$.
\begin{proof}
By Lemma \ref{relax}, $M\ba x/y$ is obtained by relaxing the circuit-hyperplane $B-y$ in $N_1/y$.  Consider $(M\ba x/y)\ba e$ for $e\notin B-y$, and assume $(M\ba x/y)\ba e$ is binary.  Then, as $(M\ba x/y)\ba e$ is a relaxation of a binary matroid, we know $(M\ba x/y)\ba e \cong U_{t-1,t}\oplus_2 U_{1,v}$ for some $t,v\geq 1$.  Now $E(M\ba x/y\ba e)$ has a partition $(S, P)$ where $S$ is the relaxed set $B-y$ and $P$ is its complement.  Then $S\cup y$ is the relaxed set of $M\ba x$ and $P\cup e$ is the relaxed set of $(M\ba x)^*$.  As $r_{M\ba x}(P\cup y)=2$, we know $|P|\leq 2$, else the matroid $N_1$ would be non-binary, and, by duality, $|S|\leq 2$.  However, $r(M\ba x)=|S\cup y|$ and $r^*(M\ba x)=|P\cup e|$.  Thus, as $M\ba x$ has a $\mathcal{W}^k$-minor for some $k\geq 3$, we know $|S\cup y|=|P\cup e|=3$.  Therefore, $M\ba x \cong \mathcal{W}^3$.  The only $3$-connected single-element extension of $\mathcal{W}^3$ that does not contain an excluded minor is $F^-_7$, depicted in Figure \ref{ext}. But $F^-_7$ is a relaxation of the Fano plane, which is binary, giving us a contradiction.  Therefore we may assume $(M\ba x/y)\ba e$ is non-binary for all $e\notin B-y$.  Then, for all such $e$, the matroid $M/y\ba e$ is also non-binary.  By Lemma \ref{relax}, for every $e\in B'$ the matroid $M/y\ba e$ is binary.  Thus if $e$ is not in $B-y$, then it is not in $B'$.  Therefore, $B'\subseteq B-y$.  As $B-y$ and $B'$ are both bases for $M/y$, the result holds. 
\end{proof}
\label{ccthyp}
\end{sublemma}

We know $B\cup e$ is a circuit of $M$ for all $e \in$ $E(M)-(B\cup x)$, and that $(B-y)\cup e$ is a circuit of $M/y$ for all $e \in$ $E(M)-B$.  Thus, since $B\cup x$ is not a circuit of $M$, we see $(B-y)\cup x$ is a circuit of $M$.  As $M\ba x$ has a $\mathcal{W}^3$-minor, but no $Q_6$-, $P_6$, or $U_{3,6}$-minor, Lemma \ref{again} implies that $M\ba x$ has a $\mathcal{W}^3$-minor using $y$.  By Lemma \ref{bases}, as $y$ is in $B$, it follows that this $\mathcal{W}^3$-minor has $y$ as a rim element.  Hence by adding $x$ back, $M$ has a single-element extension of $\mathcal{W}^3$ as a minor.  Let $\{b_1,b_2,y\}$ be the set of rim elements in $\mathcal{W}^3$.  There are only two single-element extensions of $\mathcal{W}^3$ that do not contain excluded minors, and they are $F^-_7$ and a parallel extension of a spoke element of $\mathcal{W}^3$ (see Figure \ref{ext}).  But, in each of them, the set $\{b_1,b_2,x\}$ should be a circuit because $(B-y)\cup x$ is a circuit of $M$ and the only elements of $M\ba x$ that can be contracted to produce the $\mathcal{W}^3$-minor must belong to $B$.  This contradiction completes the proof.
\end{proof}
\label{3conn}
\end{lemma}

\begin{proof}[Proofs of Theorems \ref{nb} and \ref{nb2}.]
These follow immediately by combining \\ Lemmas \ref{discon}, \ref{conn}, and \ref{3conn}. 
\end{proof}

\section{The Complexity of $\mathcal{D}$}
\label{compD}

\begin{figure}
 \begin{center}
  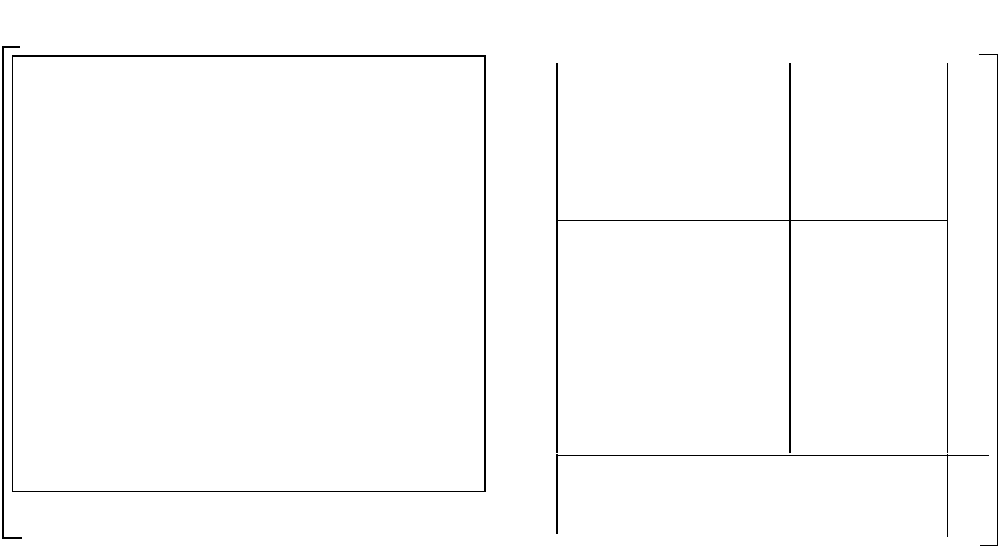
 \end{center}
\caption{The matrix $Z$.}
\label{PG}
\end{figure}

Jim Geelen asked (private communication) whether members of $\mathcal{D}$ could contain arbitrarily large projective geometries.  In this section, we observe that they can.  Note that all sums in this section are modulo two.  Let $A$ be a $k\times(2^k-1)$ matrix representing the rank-$k$ binary projective geometry $PG(k-1,2)$, where $k$ is odd.  Let $n=2^k+k+1$, let $t=2^k+k-1$, and consider the rank-$n$ binary matrix $Z$ in Figure \ref{PG}.  The entries $\alpha_i$ and $\beta_j$ are defined the next paragraph.

Let $z_{sc}$ denote the entry in row $s$ and column $c$ of $Z$.  Let $\alpha_i$=$\Sigma_{s=1}^{n-2} z_{s(n+1+i)}$, for $1\leq i\leq t$. Let $\beta_j$=$\Sigma_{c=n+2}^{2n-1} z_{jc}$, for $1\leq j\leq t$, and let $\gamma=1+\Sigma_{j=1}^t \beta_j$. Let $Z'$ be the submatrix of $Z$ whose columns are labeled by $n+1$, $n+2$, $\ldots$, $2n$.  Then each column in $Z'$ is contained in the hyperplane of $PG(n-1,2)$ consisting of those vectors whose coordinates sum to zero. Moreover, no other column of $Z$ is in this hyperplane. The definitions ensure that all the rows of $Z'$, except possibly row $n-1$, sum to zero. To see that row $n-1$ also sums to zero, note that $\Sigma_{j=1}^t \beta_j$=$\Sigma_{i=1}^t \alpha_i$ since both of these sums count the number of non-zero entries in the same submatrix.  We know that $\Sigma_{c=n+2}^{2n-1} z_{(n-1)c}$=$1+\Sigma_{i=1}^t (\alpha_i+1)$=$1+t+\Sigma_{i=1}^t \alpha_i$.  As $t$ is even, $1+t+\Sigma_{i=1}^t \alpha_i$=$1+\Sigma_{i=1}^t \alpha_i$=$1+\Sigma_{j=1}^t \beta_j=\gamma$. Thus $\{n+1,n+2,\ldots,2n\}$ is a circuit-hyperplane of $M[Z]$ and it is easy to see that its complement is as well. By relaxing both these circuit-hyperplanes, we get a member of $\mathcal{D}$ that contains a $PG(k-1,2)$-minor.

\section*{Acknowledgments}
The authors thank the anonymous referees whose many helpful suggestions significantly contributed to shaping the final version of this paper.


\begin{thebibliography}{99}

\bibitem{bixby} R.E. Bixby, $l$-matrices and a characterization of non-binary matroids, Discrete Math. 8 (1974) 139-145.

\bibitem{cun} W.H. Cunningham, A combinatorial decomposition theory, Ph.D. thesis, University of Waterloo (1973).

\bibitem{hine} N. Hine and J. Oxley, When excluding one matroid prevents infinite antichains, Adv. in Appl. Math 45 (2010) 74-76.

\bibitem{kahn} J. Kahn, A problem of P. Seymour on nonbinary matroids, Combinatorica 5 (1985) 319-323.

\bibitem{MWZ} D. Mayhew, G. Whittle, and S. van Zwam, Stability, fragility, and Rota's Conjecture, J. Combin. Theory Ser. B 102 (2012) 760-783.

\bibitem{round} J.G. Oxley, On the intersections of circuits and cocircuits in matroids, Combinatorica 4 (1984) 187-195.

\bibitem{min3} J.G. Oxley, A characterization of certain excluded-minor classes of matroids, European J. Combin. 10 (1989) 275-279.

\bibitem{list} J.G. Oxley, A characterization of a class of non-binary matroids, J. Combin. Theory Ser. B 49 (1990) 181-189.



\bibitem{book}  J. Oxley, {\it Matroid Theory}, Second edition, Oxford University Press, New York, 2011.

\bibitem{whit} J. Oxley and G. Whittle, On weak maps of ternary matroids, European J. Combin. 19 (1998) 377-389.

\bibitem{seymour} P.D. Seymour, Decomposition of regular matroids, J. Combin. Theory Ser. B 28 (1980) 305-359.

\bibitem{2round} P.D. Seymour, Minors of $3$-connected matroids, European J. Combin. 6 (1985) 375-382.



\end{thebibliography}
\end{document}